\documentclass[11pt]{article}
\usepackage{amssymb,float,amsthm,hyperref,xcolor,amsmath,tipa}
\hypersetup{
    colorlinks,
    citecolor=red,
    filecolor=blue,
    linkcolor=blue,
    urlcolor=blue
}
\usepackage[noabbrev,capitalize,nameinlink]{cleveref}
\newtheorem{theorem}{Theorem}[section]
\newtheorem{lemma}[theorem]{Lemma}

\theoremstyle{definition}

\newtheorem{question}[theorem]{Question}

\numberwithin{equation}{section}
\usepackage[margin=1in]{geometry}
\counterwithin{figure}{section}

\begin{document}
\title{On Rado numbers for equations with unit fractions}

\author{Collier Gaiser
\\
\\
\small Department of Mathematics \\
\small University of Denver\\
\small Denver, CO 80208\\
\small \tt collier.gaiser@du.edu
}

\date{}

\maketitle
\begin{abstract}
Let $f_r(k)$ be the smallest positive integer $n$ such that every $r$-coloring of $\{1,2,\ldots,n\}$ has a monochromatic solution to the nonlinear equation \[1/x_1+\cdots+1/x_k=1/y,\] where $x_1,\ldots,x_k$ are not necessarily distinct. Brown and R\"{o}dl [Bull. Aust. Math. Soc. 43(1991): 387-392] proved that $f_2(k)=O(k^6)$. In this paper, we prove that $f_2(k)=O(k^3)$. The main ingredient in our proof is a finite set $A\subseteq\mathbb{N}$ such that every $2$-coloring of $A$ has a monochromatic solution to the linear equation $x_1+\cdots+x_k=y$ and the least common multiple of $A$ is sufficiently small. This approach can also be used to study $f_r(k)$ with $r>2$. For example, a recent result of Boza, Mar\'{i}n, Revuelta, and Sanz [Discrete Appl. Math. 263(2019): 59-68] implies that $f_3(k)=O(k^{43})$.
\end{abstract}


\section{Introduction}
Let $k,r\in\mathbb{N}:=\{1,2,\ldots\}$ and $A\subseteq\mathbb{N}$. A function $\Delta:A\to\{0,1,\ldots,r-1\}$ is called an $r$-coloring of $A$, and $(x_1,\ldots,x_k,y)=(a_1,\ldots,a_k,b)$ is called a monochromatic solution to an equation $G(x_1,\ldots,x_k,y)=0$ if $G(a_1,\ldots,a_k,b)=0$ and $\Delta(a_1)=\cdots=\Delta(a_k)=\Delta(b)$.

Schur \cite{Schur1916} in 1916 proved that every finite coloring of $\mathbb{N}$ has a monochromatic solution to $x_1+x_2=y$. In 1933, Rado \cite{Rado1933} generalized this result and, among other things, proved that for all positive integers $k\geq2$, every finite coloring of $\mathbb{N}$ has a monochromatic solution to the linear equation
\begin{equation}
\sum_{i=1}^kx_i=y.
\end{equation}
By the compactness principle (see, for example, \cite{GrahamButler2015}), for all positive integers $r,k\geq2$, there exists $N\in\mathbb{N}$ such that, for all $n\geq N$, every $r$-coloring of $\{1,2,\ldots,n\}$ has a monochromatic solution to $x_1+\cdots+x_k=y$. Let $R_r(k)$ be the smallest positive integer $n$ such that every $r$-coloring of $\{1,2,\ldots,n\}$ has a monochromatic solution to $x_1+\cdots+x_k=y$, where $x_1,\ldots,x_k$ are not necessarily distinct. Beutelspacher and Brestovansky \cite{BeutelspacherBrestovansky1982} proved that $R_2(k)=k^2+k-1$ and, recently, Boza, Mar\'{i}n, Revuelta, and Sanz \cite{BMRS2019} proved that $R_3(k)=k^3+2k^2-2$.

In 1991, Brown and R\"{o}dl \cite{BrownRodl1991}, and Lefmann \cite{Lefmann1991} extended Rado's result to some nonlinear homogeneous equations. One of their results is that for all positive integers $k\geq2$, every finite coloring of $\mathbb{N}$ has a monochromatic solution to the nonlinear equation
\begin{equation}
\sum_{i=1}^k\frac{1}{x_i}=\frac{1}{y}.
\end{equation}
Let $f_r(k)$ be the smallest positive integer $n$ such that every $r$-coloring of $\{1,2,\ldots,n\}$ has a monochromatic solution to $1/x_1+\cdots+1/x_k=1/y$, where $x_1,\ldots,x_k$ are not necessarily distinct. Tejaswi and Thangdurai \cite{TejaswiThangadurai2002} proved some recursive lower bounds for $f_r(2)$ and that, for $r\geq3$, $f_r(2)\geq2^{r-3}\cdot192$. Recently, Myers and Parrish \cite{MyersParrish2018} computed that $f_2(2)=60$, $f_2(3)=40$, $f_2(4)=48$, $f_2(5)=39$, $f_3(2)=3276$, and $f_4(2)>87,000$. For general $k$, the only known result is the upper bound $f_2(k)=O(k^6)$ proved by Brown and R\"{o}dl \cite{BrownRodl1991} in 1991.
\begin{theorem}[Brown and R\"{o}dl \cite{BrownRodl1991}]\label{Theorem:BrownRodlUpper}
For all positive integers $k\geq2$,
\[
f_2(k)\leq k^2(k^2-k+1)(k^2+k-1).
\]
\end{theorem}
In this paper, we make significant improvement on \cref{Theorem:BrownRodlUpper} by showing that $f_2(k)=O(k^3)$.
\begin{theorem}\label{Theorem:NewUpper2Colors}
For all positive integers $k\geq2$,
\[
f_2(k)\leq6k(k+1)(k+2).
\]
\end{theorem}
The constant factor $6$ in \cref{Theorem:NewUpper2Colors} can be reduced to $2$ when $k\geq4$ is even or not divisible by 3. Our proof of \cref{Theorem:NewUpper2Colors} uses a variant of a theorem by Brown and R\"{o}dl \cite{BrownRodl1991}. To illustrate our idea, we first state a quantitative version of their theorem.
\begin{theorem}[Brown and R\"{o}dl \cite{BrownRodl1991}]\label{Theorem:BrownRodlLinearToNonlinear}
Let $r,k,T\in\mathbb{N}$ and $G(x_1,\ldots,x_k,y)=0$ a system of homogeneous equations such that every $r$-coloring of $\{1,2,\ldots,T\}$ has a monochromatic solution to $G(x_1,\ldots,x_k,y)=0$. Let $S$ be the least common multiple of $\{1,2,\ldots,T\}$. Then every $r$-coloring of $\{1,2,\ldots,S\}$ has a monochromatic solution to $G(1/x_1,\ldots,1/x_k,1/y)=0$.
\end{theorem}
By \cref{Theorem:BrownRodlLinearToNonlinear}, since $R_2(k)=k^2+k-1$, we have $f_2(k)\leq\text{lcm}\{1,2,\ldots,k^2+k-1\}$. It is well known that $\text{lcm}\{1,2,\ldots,k^2+k-1\}=\exp((1+o(1))k^2)$ (see, for example, Chapter 8 of \cite{Nathanson2000}). So we have $f_2(k)\leq \exp((1+o(1))k^2)$ which does not help us improve the upper bound for $f_2(k)$ in \cref{Theorem:BrownRodlUpper}. Our key observation is that the discrete interval $\{1,2,\ldots,T\}$ in \cref{Theorem:BrownRodlLinearToNonlinear} can be replaced with a finite subset of $\mathbb{N}$ whose least common multiple is smaller. Hence, to obtain a better upper bound for $f_2(k)$, it suffices to find a finite set $A\subseteq\mathbb{N}$ such that every 2-coloring of $A$ has a monochromatic solution to the linear equation $x_1+\cdots+x_k=y$ and the least common multiple of the integers in $A$ is small.

This paper is organized as follows. In \cref{Section:GeneralTheorem}, we state and prove a variant of \cref{Theorem:BrownRodlLinearToNonlinear}. In \cref{Section:ProofMainTheorem}, we prove \cref{Theorem:NewUpper2Colors}. A generalization of \cref{Theorem:NewUpper2Colors} is shown in \cref{Section:Generalizations}. Finally, we prove a polynomial upper bound for $f_3(k)$ and a lower bound for $f_r(k)$ in \cref{Section:OtherBounds}.

\subsection{Asymptotic notation}
We use standard asymptotic notation throughout this paper. For functions $f(k)$ and $g(k)$, $f(k)=O(g(k))$ if there exist constants $K$ and $C$ such that $|f(k)|\leq C|g(k)|$ for all $k\geq K$; $f(k)=\Omega(g(k))$ if there exist constants $K'$ and $c$ such that $|f(k)|\geq c|g(k)|$ for all $k\geq K'$; $f(k)=\Theta(g(k))$ if $f(k)=O(g(k))$ and $f(k)=\Omega(g(k))$; and $f(k)=o(g(k))$ if $\lim_{k\to\infty}f(k)/g(k)=0$.

\section{A variant of \cref{Theorem:BrownRodlLinearToNonlinear}}\label{Section:GeneralTheorem}
We state and prove a variant of \cref{Theorem:BrownRodlLinearToNonlinear}. This will be used in the next section for the proof of \cref{Theorem:NewUpper2Colors}.
\begin{theorem}\label{Theorem:LinearNonlinearVariant}
Let $r,k\geq2$ be integers, $A$ a finite subset of $\mathbb{N}$, $L$ the least common multiple of the integers in $A$, and $G(x_1,\ldots,x_k,y)=0$ a system of homogeneous equations. If every $r$-coloring of $A$ has a monochromatic solution to $G(x_1,\ldots,x_k,y)=0$, then every $r$-coloring of $\{1,2,\ldots,L\}$ has a monochromatic solution to $G(1/x_1,\ldots,1/x_k,1/y)=0$.
\end{theorem}
\begin{proof}
Let $r,k\geq2$ be integers, $A$ a finite subset of $\mathbb{N}$, $L$ the least common multiple of the integers in $A$, and $G(x_1,\ldots,x_k,y)=0$ a system of homogeneous equations. Suppose that every $r$-coloring of $A$ has a monochromatic solution to $G(x_1,\ldots,x_k,y)=0$. Let \[\Delta:\{1,2,\ldots,L\}\to\{0,1,\ldots,r-1\}\] be an $r$-coloring. We define \[\overline{\Delta}:A\to\{0,1,\ldots,r-1\}\] where $\overline{\Delta}(x)=\Delta(L/x)$ for all $x\in A$. Since $L$ is the least common multiple of the integers in $A$, $\overline{\Delta}$ is a well-defined function and hence an $r$-coloring of $A$. 

By assumption, there exist $a_1,\ldots,a_k,b\in A$ such that $\overline{\Delta}(a_1)=\cdots=\overline{\Delta}(a_k)=\overline{\Delta}(b)$ and $G(a_1,\ldots,a_k,b)=0$. So by our construction, $\Delta(L/a_1)=\cdots=\Delta(L/a_k)=\Delta(L/b)$. Since $G(x_1,\ldots,x_k,y)=0$ is homogeneous, we have $G(a_1/L,\ldots,a_k/L,b/L)=0$. This can be rewritten as $G(1/(L/a_1),\ldots,1/(L/a_k),1/(L/b))=0$. So $(x_1,\ldots,x_k,y)=(L/a_1,\ldots,L/a_k,L/b)$ is a monochromatic solution to $G(1/x_1,\ldots,1/x_k,1/y)=0$.
\end{proof}
\section{Proof of \cref{Theorem:NewUpper2Colors}}\label{Section:ProofMainTheorem}
We start by finding a finite subset $A$ of $\mathbb{N}$ such that every $2$-coloring of $A$ has a monochromatic solution to $x_1+\cdots+x_k=y$ and the least common multiple of the integers in $A$ is smaller than $k^2(k^2-k+1)(k^2+k-1)$ as in \cref{Theorem:BrownRodlUpper}.

\begin{lemma}\label{Lemma:KeyUpperBound}
For all integers $k\geq4$, every 2-coloring of 
\[
\{1,2,3,k+1,k+2,2(k+1),2(k+2),3k,k(k+1),k(k+2),2k(k+1),2k(k+2)\}
\]
has a monochromatic solution to $x_1+\cdots+x_k=y$.
\end{lemma}
\begin{proof}
Let $k\geq4$ be an integer. Suppose, for a contradiction, that \[\Delta:\{1,2,3,k+1,k+2,2(k+1),2(k+2),3k,k(k+1),k(k+2),2k(k+1),2k(k+2)\}\to\{R,B\}\] is a 2-coloring without a monochromatic solution to $x_1+\cdots+x_k=y$. WLOG, we assume that $\Delta(1)=R$. We have three cases depending on $\Delta(2)$ and $\Delta(3)$.

Case 1: $\Delta(2)=R$. Since 
\[
1+\cdots+1+2=(k-1)\cdot1+2=k+1
\] and \[1+\cdots+1+2+2=(k-2)\cdot1+2+2=k+2,\] we have $\Delta(k+1)=\Delta(k+2)=B$. Since \[(k+1)+\cdots+(k+1)=k\cdot(k+1)=k(k+1)\] and \[(k+2)+\cdots+(k+2)=k\cdot(k+2)=k(k+2),\] we have $\Delta(k(k+1))=\Delta(k(k+2))=R$. Now we have $\Delta(1)=\Delta(2)=\Delta(k(k+1))=\Delta(k(k+2))=R$. Since
\[
k(k+1)+1+\cdots+1+2=k(k+1)+(k-2)\cdot1+2=k(k+2),
\]
we have a monochromatic solution which is a contradiction.

Case 2: $\Delta(2)=B$ and $\Delta(3)=R$. Since \[3+\cdots+3=k\cdot3=3k\] and \[1+\cdots+1+3=(k-1)\cdot1+3=k+2,\] we have $\Delta(3k)=\Delta(k+2)=B=\Delta(2)$. Now since
\[
(k+2)+2+\cdots+2=(k+2)+(k-1)\cdot2=3k,
\]
we have a monochromatic solution which is a contradiction.

Case 3: $\Delta(2)=\Delta(3)=B$. Since \[3+3+2+\cdots+2=3+3+(k-2)\cdot2=2(k+1)\] and \[3+3+3+3+2+\cdots+2=3+3+3+3+(k-4)\cdot2=2(k+2),\] we have $\Delta(2(k+1))=\Delta(2(k+2))=R$. Since \[2(k+1)+\cdots+2(k+1)=k\cdot 2(k+1)=2k(k+1)\] and \[2(k+2)+\cdots+2(k+2)=k\cdot 2(k+2)=2k(k+2),\] we have $\Delta(2(k+1))=\Delta(2(k+2))=B$. Now we have $\Delta(2)=\Delta(3)=\Delta(2(k+1))=\Delta(2(k+2))=B$. Since
\[
2k(k+1)+3+3+2+\cdots+2=2k(k+1)+3+3+(k-3)\cdot2=2k(k+2),
\]
we have a monochromatic solution which is a contradiction.
\end{proof}

\begin{proof}[Proof of \cref{Theorem:NewUpper2Colors}]
Suppose first that $k\geq4$ is an integer. By \cref{Lemma:KeyUpperBound}, every 2-coloring of \[\{1,2,3,k+1,k+2,2(k+1),2(k+2),3k,k(k+1),k(k+2),2k(k+1),2k(k+2)\}\] has a monochromatic solution to $x_1+\cdots+x_k=y$. The least common multiple of this set of integers is at most $6k(k+1)(k+2)$. Hence, by \cref{Theorem:LinearNonlinearVariant}, $f_2(k)\leq6k(k+1)(k+2)$ for all $k\geq4$. 

By Myers and Parrish \cite{MyersParrish2018}, we also have that $f_2(2)=60<6\cdot2(2+1)(2+2)$ and $f_2(3)=40<6\cdot3(3+1)(3+2)$. Hence we have $f_2(k)\leq6k(k+1)(k+2)$ for all $k\geq2$.
\end{proof}

When $k\geq4$ is even or not divisible by $3$, we can improve the constant factor of the upper bound for $f_2(k)$ in \cref{Theorem:NewUpper2Colors} from $6$ to $2$. To see this, we first need a lemma.
\begin{lemma}\label{Lemma:UpperBoundSpecial_k}
For all even integers $k\geq4$, every 2-coloring of 
\[
\{1,2,3,k+1,k+2,2k,2(k+1),2(k+2),k(k+1),k(k+2),2k(k+1),2k(k+2)\}
\]
has a monochromatic solution to $x_1+\cdots+x_k=y$.
\end{lemma}
\begin{proof}
Let $k\geq4$ be an even integer. Suppose, for a contradiction, that \[\Delta:\{1,2,3,k+1,k+2,2k,2(k+1),2(k+2),k(k+1),k(k+2),2k(k+1),2k(k+2)\}\to\{R,B\}\] is a two-coloring without a monochromatic solution to $x_1+\cdots+x_k=y$. WLOG, we assume that $\Delta(1)=R$. We have two cases depending on $\Delta(2)$.

Case 1: $\Delta(2)=R$. The proof of this case is the same as the proof of Case 1 in \cref{Lemma:KeyUpperBound}.

Case 2: $\Delta(2)=B$. Since \[2+\cdots+2=k\cdot 2=2k,\] we have $\Delta(2k)=R$. Since \[1+\cdots+1+3+\cdots+3=(k/2)\cdot1+(k/2)\cdot3=2k,\] we have $\Delta(3)=B$. Now we have $\Delta(2)=\Delta(3)=B$ and hence the rest of the proof is the same as the proof of Case 3 in \cref{Lemma:KeyUpperBound}.
\end{proof}
\begin{theorem}\label{Theorem:NewUpperBoundSpecial_k}
If $k\geq4$ is even or not divisible by $3$, then
\[
f_2(k)\leq2k(k+1)(k+2).
\]
\end{theorem}
\begin{proof}
We first consider that $k\geq4$ is an even integer. By \cref{Lemma:UpperBoundSpecial_k}, every 2-coloring of \[\{1,2,3,k+1,k+2,2k,2(k+1),2(k+2),k(k+1),k(k+2),2k(k+1),2k(k+2)\}\] has a monochromatic solution to $x_1+\cdots+x_k=y$. Since $3$ divides $k(k+1)(k+2)$, the least common multiple of this set of integers is $2k(k+1)(k+2)$. Therefore, by \cref{Theorem:LinearNonlinearVariant}, we have $f_2(k)\leq2k(k+1)(k+2)$.

Now we suppose that $k\geq4$ and $k$ is not divisible by $3$. By \cref{Lemma:KeyUpperBound}, every 2-coloring of \[\{1,2,3,k+1,k+2,2(k+1),2(k+2),3k,k(k+1),k(k+2),2k(k+1),2k(k+2)\}\] has a monochromatic solution to $x_1+\cdots+x_k=y$. Since $3$ does not divide $k$, we have that $3$ divides $(k+1)(k+2)$ and hence the least common multiple of this set of integers is $2k(k+1)(k+2)$. Therefore, by \cref{Theorem:LinearNonlinearVariant}, we have $f_2(k)\leq2k(k+1)(k+2)$.
\end{proof}

\section{A generalization of \cref{Theorem:NewUpper2Colors}}\label{Section:Generalizations}
\cref{Theorem:NewUpper2Colors} can be used to obtain an upper bound for Rado numbers of a larger family of equations. Lefmann \cite{Lefmann1991} proved that for all integers $r,k\geq2$ and $\ell\geq1$, every $r$-coloring of $\mathbb{N}$ has a monochromatic solution to the equation
\begin{equation}
\sum_{i=1}^k\frac{1}{x_i^{1/\ell}}=\frac{1}{y^{1/\ell}}.
\end{equation}

The following result is a special case of a theorem by Lefmann \cite{Lefmann1991}. For completeness, we provide a short proof for this result.
\begin{lemma}[Lefmann \cite{Lefmann1991}]\label{Lemma:FractionalPowers}
Let $k,r\geq2$ and $\ell\geq1$ be integers, and $A$ a finite subset of $\mathbb{N}$. If every $r$-coloring of $A$ has a monochromatic solution to $1/x_1+\cdots+1/x_k=1/y$, then every $r$-coloring of $A^\ell:=\{a^\ell:a\in A\}$ has a monochromatic solution to $1/x_1^{1/\ell}+\cdots+1/x_k^{1/\ell}=1/y^{1/\ell}$.
\end{lemma}
\begin{proof}
Let $k,r\geq2$ and $\ell\geq1$ be integers, $A$ a finite subset of $\mathbb{N}$, and $A^\ell:=\{a^\ell:a\in A\}$. Suppose that every $r$-coloring of $A$ has a monochromatic solution to $1/x_1+\cdots+1/x_k=1/y$. Let
\[
\Delta:A^\ell\to\{0,1,\ldots,r-1\}
\]
be an $r$-coloring of $A^\ell$. Define
\[
\overline{\Delta}:A\to\{0,1,\ldots,r-1\}
\]
where $\overline{A}(x)=\Delta(x^\ell)$ for all $x\in A$. By the definition of $A^\ell$, $\overline{\Delta}$ a well-defined function and hence an $r$-coloring of $A$. By assumption, there exist $a_1,\ldots,a_k,b\in A$ such that $1/a_1+\cdots+1/a_k=1/b$ and $\overline{\Delta}(a_1)=\cdots=\overline{\Delta}(a_k)=\overline{\Delta}(b)$. So we have $1/(a_1^\ell)^{1/\ell}+\cdots+1/(a_k^\ell)^{1/\ell}=1/(b^\ell)^{1/\ell}$ and, by the definition of $\overline{\Delta}$, $\Delta(a_1^\ell)=\cdots=\Delta(a_k^\ell)=\Delta(b^\ell)$. Hence $(x_1,\ldots,x_k,y)=(a_1^\ell,\ldots,a_k^\ell,b^\ell)$ is a monochromatic solution to $1/x_1^{1/\ell}+\cdots+1/x_k^{1/\ell}=1/y^{1/\ell}$. 
\end{proof}

Let $f_r(k,\ell)$ be the smallest positive integer $n$ such that  every $r$-coloring of $\{1,2,\ldots,n\}$ has a monochromatic solution to $1/x_i^{1/\ell}+\cdots+1/x_k^{1/\ell}=1/y^{1/\ell}$, where $x_1,\ldots,x_k$ are not necessarily distinct. Note that we have $f_r(k,1)=f_r(k)$. The following result is a generalization of \cref{Theorem:NewUpper2Colors}.
\begin{theorem}
For all positive integers $k,r\geq2$ and $\ell\geq1$,
\[
f_2(k,\ell)\leq 6^\ell k^\ell(k+1)^\ell(k+2)^\ell.
\]
\end{theorem}
\begin{proof}
Let $k\geq2$ and $\ell\geq1$ be integers. By \cref{Theorem:NewUpper2Colors}, every 2-coloring of $\{1,2,\ldots,6k(k+1)(k+2)\}$ has a monochromatic solution to $1/x_1+\cdots+1/x_k=1/y$. Now, by \cref{Lemma:FractionalPowers}, every 2-coloring of $\{1^\ell,2^\ell,\ldots,[6k(k+1)(k+2)]^\ell\}$ has a monochromatic solution to $1/x_i^{1/\ell}+\cdots+1/x_k^{1/\ell}=1/y^{1/\ell}$. Hence we have $f_2(k,\ell)\leq[6k(k+1)(k+2)]^\ell$.
\end{proof}
Similar to \cref{Theorem:NewUpperBoundSpecial_k}, we have a slightly better upper bound for $f_2(k,\ell)$ when $k\geq4$ is even or not divisible by $3$. 
\begin{theorem}
For all positive integers $k\geq4$, $r\geq2$, and $\ell\geq1$, if $k$ is even or not divisible by $3$, then
\[
f_2(k,\ell)\leq2^\ell k^\ell(k+1)^\ell(k+2)^\ell.
\]
\end{theorem}
\begin{proof}
Let $k\geq4$, $r\geq2$, and $\ell\geq1$ be integers with $k$ even or not divisible by $3$. By \cref{Theorem:NewUpperBoundSpecial_k}, every 2-coloring of $\{1,2,\ldots,2k(k+1)(k+2)\}$ has a monochromatic solution to $1/x_1+\cdots+1/x_k=1/y$. Now, by \cref{Lemma:FractionalPowers}, every 2-coloring of $\{1^\ell,2^\ell,\ldots,[2k(k+1)(k+2)]^\ell\}$ has a monochromatic solution to $1/x_i^{1/\ell}+\cdots+1/x_k^{1/\ell}=1/y^{1/\ell}$. Hence we have $f_2(k,\ell)\leq[2k(k+1)(k+2)]^\ell$.
\end{proof}

\section{Other bounds for $f_r(k)$}\label{Section:OtherBounds}
By \cref{Theorem:LinearNonlinearVariant}, as more Ramsey-type results for linear homogeneous equations are discovered, often times related results on nonlinear homogeneous equations become direct consequences. As an example, we show that a recent result by Boza, Mar\'{i}n, Revuelta, and Sanz \cite{BMRS2019} implies a polynomial upper bound for $f_3(k)$.
\begin{lemma}[Boza, Mar\'{i}n, Revuelta, and Sanz \cite{BMRS2019}]\label{Lemma:3-Colorings}
Let $k\geq3$ and
\[
\begin{split}
\chi(k)=\{&1,2,k,k+1,k+2,2k,k^2-k+1,k^2-1,k^2,k^2+1,k^2+k-1,k^2+k,k^2+k+1,\\&2k^2-2k+1,2k^2-k,2k^2-k+1,2k^2-1,2k^2+k-2,3k^2-2k,3k^2-k-1,\\&3k^2-2,k^3,k^3+1,k^3+k-1,k^3+k,k^3+k^2-k,k^3+k^2-1,\\&k^3+k^2+k-2,k^3+2k^2-k-1,k^3+2k^2-2\}.
\end{split}
\]
Then every 3-coloring of $\chi(k)$ has a monochromatic solution to $x_1+\cdots+x_k=y$.
\end{lemma}
\begin{theorem}
$f_3(k)=O(k^{43})$.
\end{theorem}
\begin{proof}
Let $k\geq3$. The least common multiple of $\chi(k)$ as defined in \cref{Lemma:3-Colorings} is at most
\[
\begin{split}
k^3(k+1)&(k+2)(k^2-k+1)(k-1)(k^2+1)(k^2+k-1)(k^2+k+1)(2k^2-2k+1)(2k-1)\\\times &(2k^2-k+1)(2k^2-1)(2k^2+k-2)(3k-2)(3k^2-k-1)(3k^2-2)(k^3+k-1)\\&\ \ \ \ \times (k^3+k^2-1)(k^3+k^2+k-2)(k^3+2k^2-k-1)(k^3+2k^2-2)=\Theta(k^{43}).
\end{split}
\]
Hence, by \cref{Theorem:LinearNonlinearVariant}, we have $f_3(k)=O(k^{43})$.
\end{proof}
The method in this paper does not provide lower bounds for $f_r(k)$, at least not directly. Nevertheless, we note that we have the following lower bound for $f_r(k)$.
\begin{theorem}\label{Theorem:LowerBound}
For all positive integers $k,r\geq2$,
\[
f_r(k)\geq k^r.
\]
\end{theorem}
\begin{proof}
Let $k,r\geq2$ be positive integers. It suffices to show that there exists an $r$-coloring of $\{1,2,\ldots,k^r-1\}$ which does not have a monochromatic solution to $1/x_1+\cdots+1/x_k=1/y$. Consider the $r$-coloring \[\Delta:\{1,2,\ldots,k^r-1\}\to\{0,1,\ldots,r-1\}\] where $\Delta(x)=i$ if $x\in\{k^i,k^i+1,\ldots,k^{i+1}-1\}$ for some $i\in\{0,1,\ldots,r-1\}$. That is, for all $i\in\{0,1,\ldots,r-1\}$, $\Delta(\{k^i,k^i+1,\ldots,k^{i+1}-1\})=i$.

Suppose that $(x_1,\ldots,x_k,y)=(a_1,\ldots,a_k,b)$, with $b<a_1\leq\cdots\leq a_k<k^r$, is a solution to $1/x_1+\cdots+1/x_k=1/y$. Then we have
\[
\frac{1}{b}=\frac{1}{a_1}+\cdots+\frac{1}{a_k}\geq\frac{1}{a_k}+\cdots+\frac{1}{a_k}=\frac{k}{a_k}
\]
and hence $a_k\geq kb$. So $b<k^{r-1}$. It follows that $b\in\{k^j,k^j+1,\ldots,k^{j+1}-1\}$ for some $j\in\{0,1,\ldots,r-2\}$. Then $a_k\geq kb\geq kk^j=k^{j+1}$. Now by our definition, $\Delta(a_k)\geq j+1\neq j=\Delta(b)$. So $(x_1,\ldots,x_k,y)=(a_1,\ldots,a_k,b)$ is not a monochromatic solution. Therefore, $\Delta$ does not have a monochromatic solution to $1/x_1+\cdots+1/x_k=1/y$.
\end{proof}
By \cref{Theorem:LowerBound}, we have $f_2(k)=\Omega(k^2)$. Considering this and that $f_2(k)=O(k^3)$, $f_2(4)=48$, and $f_2(5)=39$, we ask the following question:
\begin{question}
Is it true that $f_2(k)=\Theta(k^2)$?
\end{question}
\section*{Declaration of competing interest}
The author declares that they have no known competing financial interests or personal relationships that could have appeared to influence the work reported in this paper.
\section*{Data availability}
No data was used for the research described in the article.
\section*{Acknowledgments}
The author would like to thank their advisor Paul Horn for continual encouragement and support. The author would also like to thank the referee for their suggestions which helped improve the presentation of the paper.

\end{document}